\documentclass[11pt]{article}
\usepackage{amsfonts}
\usepackage{mathrsfs}
\usepackage{amsmath}
\usepackage{amssymb}
\usepackage{graphicx}
\usepackage{epsfig}
\usepackage{epic}
\usepackage{cite}
\usepackage{setspace}
\usepackage{amsthm,latexsym}
\usepackage{float}

\renewcommand{\paragraph}{\roman{paragraph}}
\newtheorem{theorem}{\scshape \mdseries \bf Theorem}[section]
\newtheorem{lemma}[theorem]{\scshape \mdseries  \bf Lemma}
\newtheorem{remark}[theorem]{\scshape \mdseries  \bf Remark}

\begin{document}

\title{\sf On graphs with some normalized Laplacian eigenvalue of extremal multiplicity }
\author{ \ Fenglei Tian\thanks{Corresponding author. E-mail address: tflqsd@qfnu.edu.cn.  Supported by '' the Natural Science Foundation of Shandong Province (No. ZR2019BA016)
''. }\ ,\ \
Junqing Cai, \ Zuosong Liang, \ Xuntuan Su \\
~~\\
\noindent{\small\it \ School of Management, Qufu Normal University, Rizhao, China.}
}
\date{}
\maketitle
\noindent {\bf Abstract:} \ Let $G$ be a connected simple graph on $n$ vertices. Let $\mathcal{L}(G)$ be the normalized Laplacian matrix of $G$ and $\rho_{n-1}(G)$ be the second least eigenvalue of $\mathcal{L}(G)$. Denote by $\nu(G)$ the independence number of $G$. Recently, the paper [Characterization of graphs with some normalized Laplacian eigenvalue of multiplicity $n-3$, arXiv:1912.13227] discussed the graphs with some normalized Laplacian eigenvalue of multiplicity $n-3$. However, there is one remaining case (graphs with $\rho_{n-1}(G)\neq 1$ and $\nu(G)= 2$) not considered. In this paper, we focus on cographs and graphs with diameter 3 to investigate the graphs with some normalized Laplacian eigenvalue of multiplicity $n-3$.


\section{Introduction}

\quad { Throughout, only connected and simple graphs are considered here. Let $G=(V(G),E(G))$ be a graph with vertex set $V(G)$ and edge set $E(G)$. Let $N_G(u)$ be the set of all the neighbors of the vertex $u$. Then $d_u=|N_G(u)|$ is called the degree of $u$. For a subset $S\subset V(G)$, $S$ is called a set of twin points if $N_G(u)=N_G(v)$ for any $u, v\in S$. By $u\thicksim v$, we mean that $u$ and $v$ are adjacent. A subset $S$ of $V(G)$ is called an independent set of $G$, if the vertices of $S$ induce an empty subgraph. The cardinality of the maximum independent set of $G$ is called the independence number, denoted by $\nu(G)$. The rank of a matrix $M$ is written as $r(M)$. Let $R_{v_i}$ be the row of $M$ indexed by the vertex $v_i$.
Denote by $\mathcal{G}(n, n-3)$ the set of all $n$-vertex ($n\geq 5$) connected graphs with some normalized Laplacian eigenvalue of multiplicity $n-3$.
Let $A(G)$ and $L(G)=D(G)-A(G)$ be the adjacency matrix and the Laplacian matrix of graph $G$, respectively. Then the normalized Laplacian matrix $\mathcal{L}(G)=[l_{uv}]$ of graph $G$ is defined as $$\mathcal{L}(G)=D^{-1/2}(G)L(G)D^{-1/2}(G)=I-D^{-1/2}(G)A(G)D^{-1/2}(G),$$
where
$$l_{uv}=\begin{cases}
 \ 1, \ \ \ \ \ \ \ \ \ \ \ \ \ \text{if \ $u=v$};\\
 -1/\sqrt{d_ud_v}, \ \ \text{if $u\thicksim v$};\\
 \ 0, \ \ \ \ \ \ \ \ \ \ \ \ \ \text{otherwise}.
 \end{cases}$$
For brevity, the normalized Laplacian eigenvalues are written as $\mathcal{L}$-eigenvalues. It is well known that the least $\mathcal{L}$-eigenvalue of a connected graph is $0$ with multiplicity $1$ (see \cite{Chung}). Then let the $\mathcal{L}$-eigenvalues of a graph $G$ be
$$\rho_1(G)\geq \rho_2(G)\geq \cdots \geq \rho_{n-1}(G)> \rho_{n}(G)=0.$$

The normalized Laplacian spectrum of graphs has been studied intensively (see \cite{vanDam6,Guo,Braga,Guo1,Das,Huang2,Sun1,Sun2}), because it reveals some structural properties and some relevant dynamical aspects (such as random walk) of graphs \cite{Chung}. Recently, graphs with some eigenvalue of large multiplicity have attracted much attention (see \cite{Fernandes,Huang,Tian1}). However, there are few results on the normalized Laplacian eigenvalues. Van Dam and Omidi \cite{vanDam6} determined the graphs with some normalized Laplacian eigenvalue of multiplicity $n-1$ and $n-2$, respectively. Tian $et\ al.$ \cite{Tian} characterized two families of graphs belonging to $\mathcal{G}(n, n-3)$: graphs with $\rho_{n-1}(G)=1$ and graphs with $\rho_{n-1}(G)\neq 1$ and $\nu(G)\neq 2$, leaving the last case of $\rho_{n-1}(G)\neq 1$ and $\nu(G)= 2$ not considered. Hence, in this paper we further discuss the remaining case and obtain the following conclusion. For convenience, denote by $\mathcal{G}_1(n, n-3)$ the graphs of $\mathcal{G}(n, n-3)$ with $\rho_{n-1}(G)\neq 1$ and $\nu(G)=2$. If a graph contains no induced path $P_4$, then it is called a cograph.

\begin{theorem}\ Let $G$ be a graph of order $n\geq 5$. Then \\
(i)\ $G\in \mathcal{G}_1(n, n-3)$ with diameter 3 if and only if $G=G_1$ (see Fig. 1),\\
(ii)\ $G\in \mathcal{G}_1(n, n-3)$ and $G$ is a cograph if and only if $G=G_2$ (see Fig. 1).
\end{theorem}

\begin{figure}[htbp]
  \centering
  \setlength{\abovecaptionskip}{0cm} % 缩小caption和图像之间的距离
  \setlength{\belowcaptionskip}{0pt}
  \includegraphics[width=3.5 in]{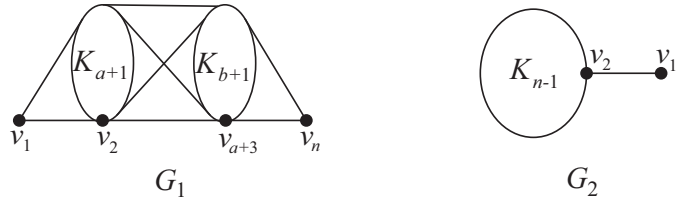}
  \caption{The graphs $G_1\ (a+b+4=n)$ and $G_2$.}
\end{figure}

Before showing the proof of Theorem 1.1, we first introduce some notations and lemmas in the next section.

\section{ Preliminaries }

For a symmetric real matrix $H$ of order $n$ whose columns and rows are indexed by $X=\{1, 2, \cdots, n\}$, let $\{X_1, X_2, \cdots, X_t\}$ be a partition of $X$.
%An $n$-dimensional column vector is called the characteristic vector of $X_i$,  if the components indexed by $X_i$ are ones and others are zeros. The $n\times t$ matrix $S$ is called the characteristic matrix of $H$, if the $i$-th column of $S$ is the characteristic vector of $X_i$.
According to the partition of $X$, we write the block form of $H$ as
$$H=
\left(
  \begin{array}{ccc}
    H_{11} & \cdots & H_{1t} \\
    \vdots & \ddots & \vdots \\
    H_{t1} & \cdots & H_{tt} \\
  \end{array}
\right),$$
where $H_{ji}$ is the transpose of $H_{ij}$.
Denote by $q_{ij}$ the average row sum of $H_{ij}$, then the matrix $Q=(q_{ij})$ is called the quotient matrix of $H$. If the row sum of $H_{ij}$ is constant, then the partition of $X$ is equitable (see \cite{Brouwer}).

\begin{lemma}\label{quotientmatrixlemma}{\rm \cite{Brouwer}} \ Suppose that $H$ is a real symmetric matrix with an equitable partition. Let $Q$ be the corresponding quotient matrix of $H$. Then, each eigenvalue of $Q$ is an eigenvalue of $H$.
\end{lemma}

\begin{lemma}\label{twinpointslemma}{\rm \cite{Huang2,Das}}\  Let $G$ be a graph with order $n$. Denote by $\{v_1,\ldots,v_p \}$ a set of twin points of $G$, then $1$ is an $\mathcal{L}$-eigenvalue of $G$ with multiplicity at least $p-1$.
\end{lemma}

\begin{lemma}\label{cliquelemma}{\rm\cite{Huang2}}\  Let $G$ be a graph with $n$ vertices. Let $K=\{v_1, \ldots, v_q \}$ be a clique in $G$ such that $N_G(v_i)-K=N_G(v_j)-K$ $(1 \leq i, j\leq q)$, then $1+\frac{1}{d_{v_i}}$ is an $\mathcal{L}$-eigenvalue of $G$ with multiplicity at least $q-1$.
\end{lemma}

\begin{lemma}\label{Tian}{\rm\cite{Tian}}\ Let $G\in\mathcal{G}(n, n-3)$ and $\theta$ be the $\mathcal{L}$-eigenvalue of $G$ with multiplicity $n-3$. If $\rho_{n-1}(G)\neq 1$, then $\theta\neq 1$.\end{lemma}

\begin{lemma}\label{path}\ Let $G\in \mathcal{G}_1(n, n-3)$ with an induced path $P_4=v_1v_2v_3v_4$. If there is a vertex $u_1$ (resp., $u_2$) such that $u_1v_2v_3v_4$ (resp., $v_1u_2v_3v_4$) is also an induced path, then $d_{v_1}=d_{u_1}$ (resp., $d_{v_2}=d_{u_2}$).
\end{lemma}

\noindent
{\bf Proof.} \ Let $\theta$ be the $\mathcal{L}$-eigenvalue of $G$ with multiplicity $n-3$,
then $r(\mathcal{L}(G)-\theta I)=3$. Since $G\in \mathcal{G}_1(n, n-3)$, it follows from Lemma \ref{Tian} that $\theta\neq 1$. Denote by $M$ the principal submatrix of $\mathcal{L}(G)-\theta I$ indexed by the vertices $\{v_1, v_2, v_3, v_4\}$, then
$$M=
\left(
  \begin{array}{cccc}
    1-\theta & \frac{-1}{\sqrt{d_{v_1}d_{v_2}}} & 0 & 0 \\
    \frac{-1}{\sqrt{d_{v_1}d_{v_2}}} & 1-\theta & \frac{-1}{\sqrt{d_{v_2}d_{v_3}}} & 0 \\
    0 & \frac{-1}{\sqrt{d_{v_2}d_{v_3}}} & 1-\theta & \frac{-1}{\sqrt{d_{v_3}d_{v_4}}} \\
    0 & 0 & \frac{-1}{\sqrt{d_{v_3}d_{v_4}}} & 1-\theta \\
  \end{array}
\right).$$
It is clear that the last three rows of $M$ are linearly independent, which yields that the rows $R_{v_2}, R_{v_3}, R_{v_4}$  of $\mathcal{L}(G)-\theta I$ are linearly independent, and then $R_{v_1}$ can be written as a linear combination of $R_{v_2}, R_{v_3}, R_{v_4}$.
Let
\begin{equation}\label{e1}
R_{v_1}=aR_{v_2}+bR_{v_3}+cR_{v_4},
\end{equation}
then
\begin{equation}\label{p1}
\begin{cases}
 \frac{-a}{\sqrt{d_{v_1}d_{v_2}}}=1-\theta\\
 a(1-\theta)-\frac{b}{\sqrt{d_{v_2}d_{v_3}}}=\frac{-1}{\sqrt{d_{v_1}d_{v_2}}}\\
 -\frac{b}{\sqrt{d_{v_3}d_{v_4}}}+c(1-\theta)=0\\
 -\frac{a}{\sqrt{d_{v_2}d_{v_3}}}+b(1-\theta)-\frac{c}{\sqrt{d_{v_3}d_{v_4}}}=0.
 \end{cases}
\end{equation}
From the first three equations of (\ref{p1}), we obtain that
\begin{equation}\label{p2}
\begin{cases}
 a=-(1-\theta)\sqrt{d_{v_1}d_{v_2}}\\
 b=-(1-\theta)^2d_{v_2}\sqrt{d_{v_1}d_{v_3}}+\sqrt{\frac{d_{v_3}}{d_{v_1}}}\\
 c=\frac{1}{(1-\theta)\sqrt{d_{v_1}d_{v_4}}}-\frac{(1-\theta)d_{v_2}\sqrt{d_{v_1}}}{\sqrt{d_{v_4}}}.
 \end{cases}
\end{equation}
Taking (\ref{p2}) into the last equation of (\ref{p1}), we deduce that
\begin{equation}\label{p3}
(1-\theta)^4d_{v_1}d_{v_2}d_{v_3}d_{v_4}-(d_{v_1}d_{v_2}+d_{v_3}d_{v_4}+d_{v_1}d_{v_4})(1-\theta)^2+1=0.
\end{equation}
Analogously, for the induced path $u_1v_2v_3v_4$ and $v_1u_2v_3v_4$, we can respectively get that
\begin{equation}\label{p4}
(1-\theta)^4d_{u_1}d_{v_2}d_{v_3}d_{v_4}-(d_{u_1}d_{v_2}+d_{v_3}d_{v_4}+d_{u_1}d_{v_4})(1-\theta)^2+1=0
\end{equation}
and
\begin{equation}\label{p5}
(1-\theta)^4d_{v_1}d_{u_2}d_{v_3}d_{v_4}-(d_{v_1}d_{u_2}+d_{v_3}d_{v_4}+d_{v_1}d_{v_4})(1-\theta)^2+1=0.
\end{equation}
Combining (\ref{p3}) and (\ref{p4}), it follows that
\begin{equation*}
(1-\theta)^2=\frac{(d_{v_2}+d_{v_4})(d_{v_1}-d_{u_1})}{d_{v_2}d_{v_3}d_{v_4}(d_{v_1}-d_{u_1})}.
\end{equation*}
If $d_{v_1}\neq d_{u_1}$, then
\begin{equation*}
(1-\theta)^2=\frac{(d_{v_2}+d_{v_4})}{d_{v_2}d_{v_3}d_{v_4}},
\end{equation*}
and thus from (\ref{p3}) we have $d_{v_3}d_{v_4}^2=0$, a contradiction. Hence, $d_{v_1}=d_{u_1}$.
Combining (\ref{p3}) and (\ref{p5}), it follows that
\begin{equation*}
(1-\theta)^2=\frac{d_{v_2}-d_{u_2}}{d_{v_3}d_{v_4}(d_{v_2}-d_{u_2})}.
\end{equation*}
If $d_{v_2}\neq d_{u_2}$, then
\begin{equation*}
(1-\theta)^2=\frac{1}{d_{v_3}d_{v_4}},
\end{equation*}
and from (\ref{p3}) we have $\frac{d_{v_1}}{d_{v_3}}=0$, a contradiction. As a result, $d_{v_2}= d_{u_2}$. The proof is completed. \hfill$\square$

\begin{lemma}\label{spectrum}\ Let $G_1$ and $G_2$ be the graphs in Fig. 1. Then their spectra (eigenvalues with multiplicity) are respectively
\begin{equation*}
\begin{cases}
\{0,\ \frac{-3+2n+\sqrt{4n-7}}{2n-4},\ \frac{-3+2n-\sqrt{4n-7}}{2n-4},\ {(\frac{n-1}{n-2})}^{n-3}\},\\
\{0,\ \frac{2n^2-5n+3+\sqrt{4n^3-23n^2+42n-23}}{2(n^2-3n+2)},\ \frac{2n^2-5n+3-\sqrt{4n^3-23n^2+42n-23}}{2(n^2-3n+2)},\ {(\frac{n-1}{n-2})}^{n-3}\}.
\end{cases}
\end{equation*}
\end{lemma}

\noindent
{\bf Proof.} \ For graph $G_1$, it is clear that $\mathcal{L}(G_1)$ has an equitable partition with respect to the vertex partition
$$V(G_1)=\{v_1\}\cup V(K_{a+1})\cup V(K_{b+1})\cup \{v_n\}.$$
Note that $d_{v_1}=a+1$, $d_{v_2}=d_{v_{a+3}}=n-2$, $d_{v_n}=b+1$ and $a+b+4=n$ in $G_1$.
Denote the quotient matrix of $\mathcal{L}(G_1)$ by $Q_1$, then
$$Q_1=\left(
         \begin{array}{cccc}
           1 & \frac{-(a+1)}{\sqrt{d_{v_1}d_{v_2}}} & 0 & 0 \\
           \frac{-1}{\sqrt{d_{v_1}d_{v_2}}} & 1-\frac{a}{d_{v_2}} & \frac{-(b+1)}{d_{v_2}d_{v_{a+3}}} & 0 \\
           0 & \frac{-(a+1)}{d_{v_2}d_{v_{a+3}}} & 1-\frac{b}{d_{v_{a+3}}} & \frac{-1}{\sqrt{d_{v_{a+3}}d_{v_n}}} \\
           0 & 0 & \frac{-(b+1)}{\sqrt{d_{v_{a+3}}d_{v_n}}} & 1 \\
         \end{array}
       \right)$$ $$=\left(
         \begin{array}{cccc}
           1 & \frac{-(a+1)}{\sqrt{(a+1)(n-2)}} & 0 & 0 \\
           \frac{-1}{\sqrt{(a+1)(n-2)}} & \frac{n-2-a}{n-2} & \frac{-(n-a-3)}{n-2} & 0 \\
           0 & \frac{-(a+1)}{n-2} & \frac{a+2}{n-2} & \frac{-1}{\sqrt{(n-a-3)(n-2)}} \\
           0 & 0 & \frac{-(n-a-3)}{\sqrt{(n-a-3)(n-2)}} & 1 \\
         \end{array}
       \right).$$
By direct calculation, the eigenvalues of $Q_1$ are
$$\{0,\ \frac{n-1}{n-2},\ \frac{-3+2n\pm\sqrt{4n-7}}{2n-4}\}.$$
From Lemma \ref{cliquelemma}, we see that $\frac{n-1}{n-2}$ is an $\mathcal{L}$-eigenvalue of $G_1$ with multiplicity at least $n-4$. It follows from the trace of $\mathcal{L}(G_1)$ and Lemma \ref{quotientmatrixlemma} that the last unknown eigenvalue is
$$n-\frac{(n-1)(n-4)}{n-2}-\frac{-3+2n-\sqrt{4n-7}}{2n-4}-\frac{-3+2n+\sqrt{4n-7}}{2n-4}=\frac{n-1}{n-2}.$$
Therefore, the multiplicity of $\frac{n-1}{n-2}$ is $n-3$.

Applying Lemma \ref{cliquelemma} to $G_2$, we obtain that $\frac{n-1}{n-2}$ is an $\mathcal{L}$-eigenvalue of $G_2$ with multiplicity at least $n-3$. Divide the vertex set $V(G_2)$ into three parts
$$V(G_2)= \{v_1\}\cup \{v_2\}\cup \{V(G_1)\setminus \{v_1,v_2\}\}.$$
Accordingly, $\mathcal{L}(G_2)$ has an equitable partition. Note that $d_{v_1}=1$, $d_{v_2}=n-1$ and the degree of each vertex of $V(G_2)\setminus \{v_1,v_2\}$ is $n-2$. Then the quotient matrix of $\mathcal{L}(G_2)$ can be written as
$$Q_2=\left(
        \begin{array}{ccc}
          1 & \frac{-1}{\sqrt{n-1}} & 0 \\
          \frac{-1}{\sqrt{n-1}} & 1 & \frac{2-n}{\sqrt{(n-1)(n-2)}} \\
          0 & \frac{-1}{\sqrt{(n-1)(n-2)}} & \frac{1}{n-2} \\
        \end{array}
      \right),$$
whose eigenvalues are
$$\{0,\ \frac{2n^2-5n+3\pm\sqrt{4n^3-23n^2+42n-23}}{2(n^2-3n+2)}\}.$$
From Lemma \ref{quotientmatrixlemma}, the eigenvalues of $Q_2$ are also the eigenvalues of $\mathcal{L}(G_2)$, which implies that the multiplicity of $\frac{n-1}{n-2}$ is $n-3$.

The proofs are completed.\hfill$\square$

\section{Main results}

\quad   Suppose that $G\in \mathcal{G}_1(n, n-3)$, $i.e.,$ $G$ contains some $\mathcal{L}$-eigenvalue of multiplicity $n-3$ with $\rho_{n-1}(G)\neq 1$ and $\nu(G)=2$, then the diameter of $G$ is not larger than 3. It is clear that the complete graph $K_n$ do not belong to $\mathcal{G}_1(n, n-3)$. Therefore, the diameter of $G$ is 2 or 3, and the proof of Theorem 1.1 is divided into Theorems 3.1 and 3.3 based on the diameter.

\begin{theorem}\label{theorem1} \ Let $G$ be a connected graph of order $n\geq 5$. Then $G\in \mathcal{G}_1(n, n-3)$ with $diam(G)=3$ if and only if $G$ is the graph $G_1$ in Fig. 1.
\end{theorem}

\noindent
{\bf Proof.} \ The sufficiency part is clear from Lemma \ref{spectrum}. In the following, we present the necessity part.
Suppose that $G\in \mathcal{G}_1(n, n-3)$ and $\theta$ is the $\mathcal{L}$-eigenvalue of multiplicity $n-3$, then $\theta\neq 1$ from Lemma \ref{Tian}. Denote by $P_4=v_1v_2v_3v_4$ a diametrical path of $G$. Assume that $U$ is a subset of $V(P_4)$ and $S_U=\{u\in V(G)\setminus V(P_4): N_G(u)\cap V(P_4)=U\}$. Since the independence number $\nu(G)=2$, we obtain that any vertex out of $V(P_4)$ must be adjacent to at least one of $V(P_4)$ and $S_{\{v_1\}}=S_{\{v_2\}}=S_{\{v_3\}}=S_{\{v_4\}}=S_{\{v_1,v_3\}}=S_{\{v_2,v_4\}}=S_{\{v_2,v_3\}}=\emptyset.$
Recalling that $diam(G)=3$, we only need to discuss the vertices of $S_{\{v_1,v_2\}}$, $S_{\{v_3,v_4\}}$, $S_{\{v_1,v_2,v_3\}}$ and $S_{\{v_2,v_3,v_4\}}$. The following claims are useful.

\vskip 2mm
\noindent
{\bf Claim 1.} \ If $S_{\{v_1,v_2\}}\neq \emptyset$, then $S_{\{v_2,v_3,v_4\}}= \emptyset$.

Suppose that $S_{\{v_2,v_3,v_4\}}\neq \emptyset$ and $v_5\in S_{\{v_1,v_2\}}$, $v_6\in S_{\{v_2,v_3,v_4\}}$. Lemma \ref{path} yields that $d_{v_1}=d_{v_5}$ and $d_{v_3}=d_{v_6}$.
Let $M_1$ be the principal submatrix of $\mathcal{L}(G)-\theta I$ indexed by $v_i\ (1\leq i\leq 6)$, then
$$M_1=\left(
      \begin{array}{cccccc}
        1-\theta & \frac{-1}{\sqrt{d_{v_1}d_{v_2}}} & 0 & 0 &  \frac{-1}{\sqrt{d_{v_1}d_{v_5}}} & 0 \\
        \frac{-1}{\sqrt{d_{v_1}d_{v_2}}} & 1-\theta & \frac{-1}{\sqrt{d_{v_2}d_{v_3}}} & 0 & \frac{-1}{\sqrt{d_{v_2}d_{v_5}}} & \frac{-1}{\sqrt{d_{v_2}d_{v_6}}} \\
        0 & \frac{-1}{\sqrt{d_{v_2}d_{v_3}}} & 1-\theta & \frac{-1}{\sqrt{d_{v_3}d_{v_4}}} & 0 &  \frac{-1}{\sqrt{d_{v_3}d_{v_6}}}\\
        0 & 0 & \frac{-1}{\sqrt{d_{v_3}d_{v_4}}} & 1-\theta & 0 & \frac{-1}{\sqrt{d_{v_4}d_{v_6}}}\\
         \frac{-1}{\sqrt{d_{v_1}d_{v_5}}} & \frac{-1}{\sqrt{d_{v_2}d_{v_5}}} & 0 & 0 & 1-\theta & *\\
        0 & \frac{-1}{\sqrt{d_{v_2}d_{v_6}}} & \frac{-1}{\sqrt{d_{v_3}d_{v_6}}} & \frac{-1}{\sqrt{d_{v_4}d_{v_6}}}& * & 1-\theta\\
      \end{array}
    \right).
$$
Clearly, the rows $R_{v_2},\ R_{v_3},\ R_{v_4}$ of $\mathcal{L}(G)-\theta I$ are linearly independent from the observation of $M_1$. Suppose the equation (\ref{e1}) still holds. Then applying (\ref{e1}) to the columns of $M_1$, we get (\ref{p1}) and
\begin{equation}\label{p6}
\begin{cases}
 -\frac{a}{\sqrt{d_{v_2}d_{v_5}}}=-\frac{1}{\sqrt{d_{v_1}d_{v_5}}}\\
 -\frac{a}{\sqrt{d_{v_2}d_{v_6}}}-\frac{b}{\sqrt{d_{v_3}d_{v_6}}}-\frac{c}{\sqrt{d_{v_4}d_{v_6}}}=0.
 \end{cases}
\end{equation}
The first equations of (\ref{p1}) and (\ref{p6}) indicate that $1-\theta=-\frac{1}{d_{v_1}}$, which together with (\ref{p3}) yields that
\begin{equation}\label{p7}
  d_{v_1}^2(d_{v_1}-d_{v_2}-d_{v_4})=d_{v_3}d_{v_4}(d_{v_1}-d_{v_2}).
\end{equation}
The fourth equation of (\ref{p1}) and the second one of (\ref{p6}) imply that $\frac{b}{\sqrt{d_{v_3}}}=\frac{b\sqrt{d_{v_3}}}{d_{v_1}}$. Further, as $b\neq 0$ (otherwise $a=c=0$ from (\ref{p1}), a contradiction), then we get
\begin{equation}\label{p8}
d_{v_1}=d_{v_3}.
\end{equation}
Bringing (\ref{p8}) into (\ref{p7}), we derive $d_{v_3}^2(d_{v_1}-d_{v_2}-d_{v_4})=d_{v_3}d_{v_4}(d_{v_1}-d_{v_2})$, that is,
\begin{equation}\label{p9}
(d_{v_3}-d_{v_4})(d_{v_1}-d_{v_2})=d_{v_3}d_{v_4}.
\end{equation}
Note that each vertex of $V(G)\setminus V(P_4)$ just belongs to $S_{\{v_1,v_2\}}$, $S_{\{v_3,v_4\}}$, $S_{\{v_1,v_2,v_3\}}$ or $S_{\{v_2,v_3,v_4\}}$. Then it is easy to know that $d_{v_2}>d_{v_1}$ and $d_{v_3}>d_{v_4}$.
As a result, the left side of (\ref{p9}) is negative, but the right side is positive, a contradiction.

\vskip 2mm
\noindent
{\bf Claim 2.} \ $S_{\{v_1,v_2\}}=S_{\{v_3,v_4\}}=\emptyset$.

By symmetry, it suffices to prove $S_{\{v_1,v_2\}}= \emptyset$. Suppose on the contrary that $S_{\{v_1,v_2\}}\neq \emptyset$ and $v_5\in S_{\{v_1,v_2\}}$. Let $M_2$ be the principal submatrix of $\mathcal{L}(G)-\theta I$ indexed by $v_i\ (1\leq i\leq 5)$, then $M_2$ is a principal submatrix of $M_1$ and the equation (\ref{p7}) still holds. From Claim 1, we see that $S_{\{v_2,v_3,v_4\}}=\emptyset$, which implies that $d_{v_2}=d_{v_1}+1$. Further, if $S_{\{v_3,v_4\}}=\emptyset$ at this moment, then $d_{v_4}=1$. Reconsidering (\ref{p7}), we obtain that $2d_{v_1}^2=d_{v_3}$, contradicting with the fact that $d_{v_1}\geq d_{v_3}$. On the other hand, suppose $S_{\{v_3,v_4\}}\neq \emptyset$ and $v_6\in S_{\{v_3,v_4\}}$, then $S_{\{v_1,v_2,v_3\}}=\emptyset$ from Claim 1 and $d_{v_3}=d_{v_4}+1$. For the subgraph induced by $\{v_1,v_2,v_3,v_4,v_6\}$, similar deduction with (\ref{p7}) leads to
\begin{equation}\label{p10}
  d_{v_4}^2(d_{v_4}-d_{v_3}-d_{v_1})=d_{v_1}d_{v_2}(d_{v_4}-d_{v_3}).
\end{equation}
Applying $d_{v_2}=d_{v_1}+1$ and $d_{v_3}=d_{v_4}+1$ to (\ref{p7}) and (\ref{p10}), we derive that $d_{v_1}^2=d_{v_4}$ and $d_{v_4}^2=d_{v_1}$, which imply that $d_{v_1}=d_{v_4}=1$, a contradiction. Therefore, $S_{\{v_1,v_2\}}= \emptyset$.

Next, we show that $G$ must be isomorphic to $G_2$ in Fig. 1.  From Claim 2, the vertices of $V(G)\setminus V(P_4)$ belong to $S_{\{v_1,v_2,v_3\}}$ or $S_{\{v_2,v_3,v_4\}}$. First, any two vertices of $S_{\{v_1,v_2,v_3\}}$
(resp., $S_{\{v_2,v_3,v_4\}}$) are adjacent. Otherwise, it is easy to see that $\nu(G)\geq 3$, contradicting with $\nu(G)=2$. Further, suppose $u\in S_{\{v_1,v_2,v_3\}}$ (resp., $S_{\{v_2,v_3,v_4\}}$), then $d_u=d_{v_2}$ (resp., $d_u=d_{v_3}$) by Lemma \ref{path}, which indicates that $u$ is adjacent to each of $S_{\{v_2,v_3,v_4\}}$ (resp., $S_{\{v_1,v_2,v_3\}}$). Hence, $G$ is isomorphic to $G_2$.

The proof is completed. \hfill$\square$

Now, we discuss the case of $diam(G)=2$. In the following, we always let $P_3=v_1v_2v_3$ be a diametrical path of $G$. Assume that $U$ is a subset of $V(P_3)$ and
$$S_U=\{u\in V(G)\setminus V(P_3): N_G(u)\cap V(P_3)=U\}.$$

\begin{lemma}\label{lemma1}\ Let $G$ be a cograph with $diam(G)=2$ and $G\in \mathcal{G}_1(n, n-3)$. The diametrical path $P_3$ and the notation $S_U$ are stated as above. Then \\
(i)\ $S_{\{v_1\}}=S_{\{v_2\}}=S_{\{v_3\}}=\emptyset$;\\
(ii)\ If $S_{\{v_1,v_2\}}\neq \emptyset$, then $S_{\{v_1,v_3\}}=S_{\{v_2,v_3\}}= S_{\{v_1,v_2,v_3\}}=\emptyset$.
\end{lemma}

\noindent
{\bf Proof.} \ Since $\nu(G)=2$, then each vertex out of $V(P_3)$ must be adjacent to $v_1$ or $v_3$. Hence, $S_{\{v_2\}}=\emptyset$.
Note that $G$ is a cograph (i.e., containing no induced path $P_4$), then $S_{\{v_1\}}=S_{\{v_3\}}=\emptyset$.
Therefore, the vertices out of $V(G)\setminus V(P_3)$ belong to $S_{\{v_1,v_2\}}$, $S_{\{v_1,v_3\}}$, $S_{\{v_2,v_3\}}$ or $S_{\{v_1,v_2,v_3\}}$. The remaining proof can be completed by the following claims.

\vskip 2mm
\noindent
{\bf Claim 1.} \ If $S_{\{v_1,v_2\}}\neq \emptyset$, then $S_{\{v_1,v_3\}}=\emptyset$. In other words, $G$ contains no induced subgraph isomorphic to $H_1$ (see Fig. 2).

Let $v_4\in S_{\{v_1,v_2\}}$. Suppose for a contradiction that $S_{\{v_1,v_3\}}\neq \emptyset$ and $v_5\in S_{\{v_1,v_3\}}$. Then $v_4\thicksim v_5$, otherwise $\{v_4,v_1,v_5,v_3\}$ induce a path $P_4$, a contradiction. Thus, the vertices $v_i\ (1\leq i\leq 5)$ induce a subgraph of $G$ isomorphic to $H_1$ in Fig. 2. Denote by $M_2$ the principal submatrix of $\mathcal{L}(G)-\theta I$ indexed by $v_i\ (1\leq i\leq 5)$, then
$$M_2=
\left(
  \begin{array}{ccccc}
    1-\theta & \frac{-1}{\sqrt{d_{v_1}d_{v_2}}} & 0 & \frac{-1}{\sqrt{d_{v_1}d_{v_4}}} & \frac{-1}{\sqrt{d_{v_1}d_{v_5}}} \\
    \frac{-1}{\sqrt{d_{v_1}d_{v_2}}} & 1-\theta & \frac{-1}{\sqrt{d_{v_2}d_{v_3}}} & \frac{-1}{\sqrt{d_{v_2}d_{v_4}}} & 0 \\
    0 & \frac{-1}{\sqrt{d_{v_2}d_{v_3}}} & 1-\theta & 0 & \frac{-1}{\sqrt{d_{v_3}d_{v_5}}} \\
    \frac{-1}{\sqrt{d_{v_1}d_{v_4}}} & \frac{-1}{\sqrt{d_{v_2}d_{v_4}}} & 0 & 1-\theta & \frac{-1}{\sqrt{d_{v_4}d_{v_5}}} \\
    \frac{-1}{\sqrt{d_{v_1}d_{v_5}}} & 0 & \frac{-1}{\sqrt{d_{v_3}d_{v_5}}} & \frac{-1}{\sqrt{d_{v_4}d_{v_5}}} & 1-\theta \\
  \end{array}
\right).
$$
Since the following minor $D$ of $M_2$ is nonzero,
$$D=\left|\begin{array}{ccc}
         \frac{-1}{\sqrt{d_{v_2}d_{v_3}}} & \frac{-1}{\sqrt{d_{v_2}d_{v_4}}} & 0 \\
         1-\theta & 0 & \frac{-1}{\sqrt{d_{v_3}d_{v_5}}} \\
         0 & 1-\theta & \frac{-1}{\sqrt{d_{v_4}d_{v_5}}} \\
       \end{array}\right|=
       -\frac{1-\theta}{\sqrt{d_{v_2}d_{v_5}}}(\frac{1}{d_{v_3}}+\frac{1}{d_{v_4}})\neq 0,$$
then the second, third and fourth rows of $M_2$ are linearly independent, and thus the rows $R_{v_2},\ R_{v_3},\ R_{v_4}$ of $\mathcal{L}(G)-\theta I$ are linearly independent. As $r(\mathcal{L}(G)-\theta I)=3$, any row of $\mathcal{L}(G)-\theta I$ can be represented as a linear combination of $R_{v_i}\ (2\leq i\leq 4)$. Let
\begin{equation}\label{p11}
R_{v_1}=aR_{v_2}+bR_{v_3}+cR_{v_4}.
\end{equation}
Applying (\ref{p11}) to the columns of $M_2$, we get
\begin{equation}\label{p12}
\begin{cases}
 -\frac{a}{\sqrt{d_{v_1}d_{v_2}}}-\frac{c}{\sqrt{d_{v_1}d_{v_4}}}=1-\theta\\
  a(1-\theta)-\frac{b}{\sqrt{d_{v_2}d_{v_3}}}-\frac{c}{\sqrt{d_{v_2}d_{v_4}}}=-\frac{1}{\sqrt{d_{v_1}d_{v_2}}}\\
 -\frac{a}{\sqrt{d_{v_2}d_{v_3}}}+b(1-\theta)=0\\
 -\frac{a}{\sqrt{d_{v_2}d_{v_4}}}+c(1-\theta)=-\frac{1}{\sqrt{d_{v_1}d_{v_4}}}\\
 -\frac{b}{\sqrt{d_{v_3}d_{v_5}}}-\frac{c}{\sqrt{d_{v_4}d_{v_5}}}=-\frac{1}{\sqrt{d_{v_1}d_{v_5}}}.
 \end{cases}
\end{equation}
The second and the fifth equations of (\ref{p12}) imply that $a=0$. Then $b=0$ from the third equation of (\ref{p12}), which indicates that $c=\sqrt{\frac{d_{v_4}}{d_{v_1}}}$ by the fifth one of (\ref{p12}). Taking the values of $a$ and $c$ into the first and the fourth ones of (\ref{p12}), we obtain that
\begin{equation}\label{p13}
1-\theta=-\frac{1}{d_{v_1}}=-\frac{1}{d_{v_4}},
\end{equation}
which yields that $d_{v_1}=d_{v_4}$. Moreover, let
\begin{equation*}
R_{v_5}=sR_{v_2}+tR_{v_3}+kR_{v_4},
\end{equation*}
and we have the following equations from the columns of $M_2$,
\begin{equation}\label{p14}
\begin{cases}
 -\frac{s}{\sqrt{d_{v_1}d_{v_2}}}-\frac{k}{\sqrt{d_{v_1}d_{v_4}}}=-\frac{1}{\sqrt{d_{v_1}d_{v_5}}}\\
  s(1-\theta)-\frac{t}{\sqrt{d_{v_2}d_{v_3}}}-\frac{k}{\sqrt{d_{v_2}d_{v_4}}}=0\\
 -\frac{s}{\sqrt{d_{v_2}d_{v_3}}}+t(1-\theta)=-\frac{1}{\sqrt{d_{v_3}d_{v_5}}}\\
 -\frac{s}{\sqrt{d_{v_2}d_{v_4}}}+k(1-\theta)=-\frac{1}{\sqrt{d_{v_4}d_{v_5}}}\\
 -\frac{t}{\sqrt{d_{v_3}d_{v_5}}}-\frac{k}{\sqrt{d_{v_4}d_{v_5}}}=1-\theta.
 \end{cases}
\end{equation}
It follows from the third and the fourth equations of (\ref{p14}) that $t\sqrt{d_{v_3}}=k\sqrt{d_{v_4}}$, which, together with (\ref{p13}) and the fifth one of (\ref{p14}), implies that
\begin{equation*}
\begin{cases}
 t=\frac{\sqrt{d_{v_3}d_{v_5}}}{d_{v_3}+d_{v_4}}\\
 k=\frac{d_{v_3}\sqrt{d_{v_5}}}{\sqrt{d_{v_4}}(d_{v_3}+d_{v_4})}.
 \end{cases}
\end{equation*}
Bringing the value of $k$ into the first one of (\ref{p14}), we have $s=\sqrt{\frac{d_{v_2}}{d_{v_5}}}-\frac{d_{v_3}\sqrt{d_{v_2}d_{v_5}}}{d_{v_4}(d_{v_3}+d_{v_4})}$. Then now the second one of (\ref{p14}) can be simplified to be
\begin{equation}\label{p15}
  d_{v_4}^2d_{v_2}+d_{v_4}^2d_{v_5}+d_{v_2}d_{v_3}d_{v_4}+d_{v_3}d_{v_4}d_{v_5}=d_{v_2}d_{v_3}d_{v_5},
\end{equation}
implying that $d_{v_2}>d_{v_4}$. Thus it follows from $d_{v_1}=d_{v_4}$ that $d_{v_2}>d_{v_1}$, which indicates that $S_{\{v_2,v_3\}}\neq \emptyset$. Let $v_6\in S_{\{v_2,v_3\}}$, then $v_5\thicksim v_6$ (otherwise $\{v_1,v_5,v_3,v_6\}$ induce $P_4$). From the symmetry of $v_4$ and $v_6$, applying similar discussion to the subgraph induced by $\{v_1,v_2,v_3,v_5,v_6\}$, one can obtain that
\begin{equation}\label{p16}
  1-\theta=-\frac{1}{d_{v_3}}=-\frac{1}{d_{v_6}}.
\end{equation}
Combining (\ref{p13}) and (\ref{p16}), $d_{v_3}=d_{v_4}$ holds. Then the equation (\ref{p15}) can be rewritten as
\begin{equation}\label{p17}
  d_{v_4}d_{v_2}+d_{v_4}d_{v_5}+d_{v_2}d_{v_3}+d_{v_3}d_{v_5}=d_{v_2}d_{v_5}.
\end{equation}
Recalling that any vertex out of $V(P_3)$ must be adjacent to $v_1$ or $v_3$, we can see $d_{v_1}+d_{v_3}>d_{v_2}$, i.e., $d_{v_4}+d_{v_3}>d_{v_2}$ as $d_{v_1}=d_{v_4}$, contradicting with (\ref{p17}).
As a result, if $S_{\{v_1,v_2\}}\neq \emptyset$, then $S_{\{v_1,v_3\}}= \emptyset$.

\begin{figure}[htbp]
  \centering
  \setlength{\abovecaptionskip}{0cm} % 缩小caption和图像之间的距离
  \setlength{\belowcaptionskip}{0pt}
  \includegraphics[width=6 in]{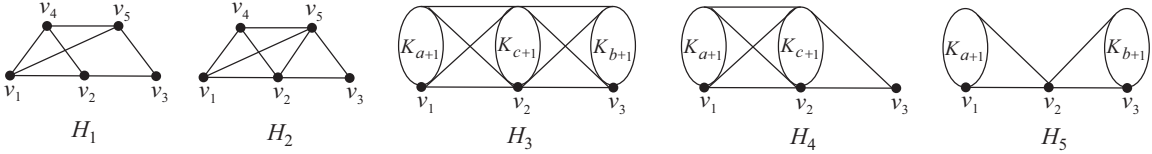}
  \caption{The graphs $H_1, H_2$, $H_3 (a+b+c+3=n)$, $H_4 (a+c+3=n)$ and $H_5 (a+b+3=n)$.}
\end{figure}

\vskip 2mm
\noindent
{\bf Claim 2.} \ If $S_{\{v_1,v_2\}}\neq \emptyset$, then $S_{\{v_1,v_2,v_3\}}=\emptyset$. In other words, $G$ contains no induced subgraph isomorphic to $H_2$ (see Fig. 2).

Suppose on the contrary that $S_{\{v_1,v_2,v_3\}}\neq \emptyset$ and $v_5\in S_{\{v_1,v_2,v_3\}}$. Also, let $v_4\in S_{\{v_1,v_2\}}$. It is easy to see that $v_4\thicksim v_5$, otherwise $\{v_4,v_1,v_5,v_3\}$ induce a path $P_4$, a contradiction. Then $G$ contains $H_2$ as an induced subgraph.
Since $S_{\{v_1,v_3\}}=\emptyset$ from Claim 1, then all other vertices maybe spread in $S_{\{v_1,v_2\}}$, $S_{\{v_1,v_2,v_3\}}$ or $S_{\{v_2,v_3\}}$.

\vskip 1.5mm
\noindent
{\sf Case 1.} \  Assume that $S_{\{v_2,v_3\}}\neq \emptyset$.

We point out that all vertices of $S_{\{v_1,v_2\}}$ (resp., $S_{\{v_2,v_3\}}$) induce a clique of $G$, otherwise one can see $\nu(G)\geq 3$, a contradiction. Further, each vertex of $S_{\{v_1,v_2\}}$ (resp., $S_{\{v_2,v_3\}}$) is adjacent to each one of $S_{\{v_1,v_2,v_3\}}$, otherwise one can easily obtain an induced path $P_4$, a contradiction.
Additionally, all vertices of $S_{\{v_1,v_2,v_3\}}$ also induce a clique. If not, let $u,w\in S_{\{v_1,v_2,v_3\}}$ and $u\nsim w$, then $\{v_1,v_3,v_4,u,w\}$ induce a subgraph isomorphic to $H_1$ in Fig. 2, contradicting with Claim 1.
Let $|S_{\{v_1,v_2\}}|=a\geq 1$, $|S_{\{v_2,v_3\}}|=b\geq 1$ and $|S_{\{v_1,v_2,v_3\}}|=c\geq 1$, then $G$ is isomorphic to $H_3\ (a+b+c+3=n)$  in Fig. 2. The remaining proof is divided into the following cases.

\vskip 1.5mm
\noindent
$\bullet$\  Suppose that $a=b$ and $a+b\leq c$ in $H_3$. Then we declare that $a=b=1$. Otherwise, $a=b\geq 2$, and then $c\geq a+b\geq 4$. From Lemma \ref{cliquelemma}, $1+\frac{1}{d_{v_1}}=1+\frac{1}{d_{v_3}}$ (resp., $1+\frac{1}{d_{v_2}}$) is an $\mathcal{L}$-eigenvalue with multiplicity at least $4$. Noting that $d_{v_1}=d_{v_3}\neq d_{v_2}$, then $G$ contains no eigenvalue of multiplicity $n-3$, a contradiction. Hence, $a=b=1$ and $c\geq 2$. By Lemma \ref{cliquelemma} again, the multiplicity of $1+\frac{1}{d_{v_1}}$ (resp., $1+\frac{1}{d_{v_2}}$) is at least 2.
As a result,
$$\theta=1+\frac{1}{d_{v_1}}=1+\frac{1}{n-3}\ \ \text{or}\ \ \theta=1+\frac{1}{d_{v_2}}=1+\frac{1}{n-1}.$$
If $\theta=1+\frac{1}{n-3}$, then the multiplicity of  $1+\frac{1}{n-1}$ is equal to 2. Then from the trace of $\mathcal{L}(G)$,
\begin{equation*}
  n = (n-3)(1+\frac{1}{n-3})+2(1+\frac{1}{n-1})
   = n+\frac{2}{n-1}>n,
\end{equation*}
a contradiction.
If $\theta=1+\frac{1}{n-1}$, then the multiplicity of $1+\frac{1}{n-3}$ is equal to 2. Then
\begin{equation*}
  n = (n-3)(1+\frac{1}{n-1})+2(1+\frac{1}{n-3})
   = n+\frac{4}{n^2-4n+3}>n,
\end{equation*}
a contradiction.

\vskip 1.5mm
\noindent
$\bullet$\  Suppose that $a=b$ and $a+b> c$ in $H_3$. Then we say that $c=1$. Otherwise, $c\geq 2$ and $a+b\geq 4$.
By Lemma \ref{cliquelemma}, the multiplicity of $1+\frac{1}{d_{v_1}}$ (resp., $1+\frac{1}{d_{v_2}}$) as an $\mathcal{L}$-eigenvalue is at least $4$ (resp., at least $2$). Thus, $\theta=1+\frac{1}{d_{v_1}}$ and the multiplicity of $1+\frac{1}{d_{v_2}}$ is equal to 2, implying that $c=2$. Note that $n=a+b+c+3>7$ in this case, then by the trace of $\mathcal{L}(G)$,
\begin{equation*}
  \begin{array}{rcl}
  n &=& (n-3)(1+\frac{1}{d_{v_1}})+2(1+\frac{1}{d_{v_2}})\\
   &=& (n-3)(1+\frac{2}{n+1})+2(1+\frac{1}{n-1})\\
   &=& n+\frac{n-7}{n+1}+\frac{2}{n-1}>n
    \end{array}
\end{equation*}
a contradiction. Thus, $c=1$ and $a=b=\frac{n-4}{2}$.
Clearly,  $\mathcal{L}(G)$ has an equitable partition according to $V(G)=\{V(K_{a+1}),V(K_{b+1}),V(K_{c+1})\}$. Let the quotient matrix of $\mathcal{L}(G)$ be $Q_1$, then from $d_{v_1}=d_{v_3}=\frac{n}{2}$ and $d_{v_2}=n-1$
$$Q_1=\left(
      \begin{array}{ccc}
        1-\frac{n-4}{n} & \frac{-2\sqrt{2}}{\sqrt{n(n-1)}} & 0 \\
        \frac{2-n}{\sqrt{2n(n-1)}} & 1-\frac{1}{n-1} & \frac{2-n}{\sqrt{2n(n-1)}} \\
        0 & \frac{-2\sqrt{2}}{\sqrt{n(n-1)}} & 1-\frac{n-4}{n} \\
      \end{array}
    \right).
$$
By direct calculation, the eigenvalues of $Q_1$ are $\{0, \frac{4}{n},\frac{n^2+2n-4}{n^2-n}\}$.
By Lemma \ref{cliquelemma}, $1+\frac{1}{d_{v_1}}=\frac{n+2}{n}$ and $1+\frac{1}{d_{v_2}}=\frac{n}{n-1}$ are two distinct $\mathcal{L}$-eigenvalues of $G$. Then from Lemma \ref{quotientmatrixlemma}, we see that $G$ has 5 distinct $\mathcal{L}$-eigenvalues, contradicting with $G\in \mathcal{G}_1(n, n-3)$.

\vskip 1.5mm
\noindent
$\bullet$\  Suppose that $a\neq b$ and $a>b$ without loss of generality. From $b\geq 1$ and $c\geq 1$, then $a\geq 2$ and $1+\frac{1}{d_{v_1}}$ ( with multiplicity at least 2), $1+\frac{1}{d_{v_2}}$ and $1+\frac{1}{d_{v_3}}$ are three distinct $\mathcal{L}$-eigenvalues of $G$ by Lemma \ref{cliquelemma}. Thus, one can derive that $\theta=1+\frac{1}{d_{v_1}}$, which yields that $c=b=1$. Therefore, the trace of $\mathcal{L}(G)$ is
\begin{equation*}
  \begin{array}{rcl}
  n &=& (n-3)(1+\frac{1}{d_{v_1}})+(1+\frac{1}{d_{v_2}})+(1+\frac{1}{d_{v_3}})\\
    &=& (n-3)(1+\frac{1}{n-3})+(1+\frac{1}{n-1})+(1+\frac{1}{3})\\
    &=& n+\frac{1}{n-1}+\frac{1}{3}>n,
  \end{array}
\end{equation*}
a contradiction. From above three subcases for $S_{\{v_2,v_3\}}\neq \emptyset$, we can always obtain contradictions.

\vskip 1.5mm
\noindent
{\sf Case 2.} \  Let $S_{\{v_2,v_3\}}= \emptyset$.

In this case, there are only $S_{\{v_1,v_2\}}$ and $S_{\{v_1,v_2,v_3\}}$ nonempty.
Similar as Case 1, all vertices of $S_{\{v_1,v_2\}}$ (resp., $S_{\{v_1,v_2,v_3\}}$) induce a clique of $G$ and each vertex of $S_{\{v_1,v_2\}}$ is adjacent to each one of $S_{\{v_1,v_2,v_3\}}$. Let $|S_{\{v_1,v_2\}}|=a\geq 1$ and $|S_{\{v_1,v_2,v_3\}}|=c\geq 1$, then $G$ is isomorphic to $H_4$ in Fig. 2.

\vskip 1.5mm
\noindent
$\bullet$\  Suppose that $a\geq c$, then we claim that $c\leq 2$ by lemma \ref{cliquelemma}. If $c=2$, then $1+\frac{1}{d_{v_1}}$ and $1+\frac{1}{d_{v_2}}$ are two distinct $\mathcal{L}$-eigenvalues of $G$ with multiplicity at least 2 by Lemma \ref{cliquelemma}. Thus, $\theta=1+\frac{1}{d_{v_1}}$ or $\theta=1+\frac{1}{d_{v_2}}$.
If $\theta=1+\frac{1}{d_{v_1}}$, then the multiplicity of $1+\frac{1}{d_{v_2}}$ is 2 and from the trace of $\mathcal{L}(G)$,
\begin{equation*}
  \begin{array}{rcl}
  n&=&(n-3)(1+\frac{1}{d_{v_1}})+2(1+\frac{1}{d_{v_2}})\\
  &=&(n-3)(1+\frac{1}{n-2})+2(1+\frac{1}{n-1})\\
  &=&n+\frac{n-3}{(n-1)(n-2)}>n,
    \end{array}
\end{equation*}
a contradiction. If $\theta=1+\frac{1}{d_{v_2}}$, then the multiplicity of $1+\frac{1}{d_{v_1}}$ is 2, which yields that
$a=2$ (as $a\geq c=2$). Then $G$ is a graph isomorphic to $H_4$ with order $7$ and by direct check $G\notin \mathcal{G}_1(n, n-3)$. If $c=1$, then $a=n-4$. According to $V(G)=\{V(K_{a+1}), V(K_{c+1}), v_3\}$, there is an equitable partition for $\mathcal{L}(G)$. Let $Q_2$ be the corresponding quotient matrix of $\mathcal{L}(G)$, then by $d_{v_1}=n-2$, $d_{v_2}=n-1$ and $d_{v_3}=2$,
$$Q_2=
\left(
  \begin{array}{ccc}
    1-\frac{n-4}{n-2} & \frac{-2}{\sqrt{(n-2)(n-1)}} & 0 \\
    \frac{3-n}{\sqrt{(n-2)(n-1)}} & 1-\frac{1}{n-1} & \frac{-1}{\sqrt{2(n-1)}} \\
    0 & \frac{-2}{\sqrt{2(n-1)}} & 1 \\
  \end{array}
\right).
$$
By calculation, the eigenvalues of $Q_2$ are
$$\{0,\ \frac{2n^2-5n+4\pm\sqrt{4n^3-27n^2+56n-32}}{2(n^2-3n+2)}\},$$
which are also the eigenvalues of $\mathcal{L}(G)$ from Lemma \ref{quotientmatrixlemma}.
Furthermore, $1+\frac{1}{d_{v_1}}=\frac{n-1}{n-2}$ and $1+\frac{1}{d_{v_2}}=\frac{n}{n-1}$ are two distinct $\mathcal{L}$-eigenvalues of $G$ by Lemma \ref{cliquelemma}. Thus $G$ has 5 distinct $\mathcal{L}$-eigenvalues, contradicting with $G\in \mathcal{G}_1(n, n-3)$.

\vskip 1.5mm
\noindent
$\bullet$\  Suppose that $a< c$, then $a\leq 2$. Otherwise, $a\geq 3$ and the multiplicities of $1+\frac{1}{d_{v_1}}$ and $1+\frac{1}{d_{v_2}}$ as two distinct $\mathcal{L}$-eigenvalues of $G$ are at least 3 by Lemma \ref{cliquelemma}, contradicting with $G\in \mathcal{G}_1(n, n-3)$. First, assume that $a=2$, then $c\geq 3$ and $\theta=1+\frac{1}{d_{v_2}}$ clearly, which implies that the multiplicity of $1+\frac{1}{d_{v_1}}$ must be 2. Thus, from the trace of $\mathcal{L}(G)$, we obtain
\begin{equation*}
  \begin{array}{rcl}
  n&=&2(1+\frac{1}{d_{v_1}})+(n-3)(1+\frac{1}{d_{v_2}})\\
  &=&2(1+\frac{1}{n-2})+(n-3)(1+\frac{1}{n-1})\\
  &=&n+\frac{2}{(n-1)(n-2)}>n,
    \end{array}
\end{equation*}
a contradiction. Now, assume that $a=1$, then $c=n-4$. It is clear that $\mathcal{L}(G)$ contains an equitable partition with respect to $V(G)=\{V(K_{a+1}), V(K_{c+1}), v_3\}$. Note that $d_{v_1}=n-2$, $d_{v_2}=n-1$ and $d_{v_3}=n-3$, then the corresponding quotient matrix of $\mathcal{L}(G)$ is
$$Q_3=
\left(
  \begin{array}{ccc}
    1-\frac{1}{n-2} & \frac{3-n}{\sqrt{(n-2)(n-1)}} & 0 \\
    \frac{-2}{\sqrt{(n-2)(n-1)}} & 1-\frac{n-4}{n-1} & \frac{-1}{\sqrt{(n-1)(n-3)}} \\
    0 & \frac{3-n}{\sqrt{(n-1)(n-3)}} & 1 \\
  \end{array}
\right).$$
Further, by calculating, the eigenvalues of $Q_3$ are
$$\{0,\ \frac{2n^2-4n-1\pm\sqrt{8n^2-32n+33}}{2(n^2-3n+2)}\},$$
which are also the eigenvalues of $\mathcal{L}(G)$ by Lemma \ref{quotientmatrixlemma}.
Recalling that $1+\frac{1}{d_{v_1}}=\frac{n-1}{n-2}$ and $1+\frac{1}{d_{v_2}}=\frac{n}{n-1}$ are two distinct $\mathcal{L}$-eigenvalues of $G$, then we see that $G$ has 5 distinct $\mathcal{L}$-eigenvalues, contradicting with $G\in \mathcal{G}_1(n, n-3)$. From above two subcases for $S_{\{v_2,v_3\}}= \emptyset$, we can also obtain contradictions.

Consequently, we conclude that if $S_{\{v_1,v_2\}}\neq \emptyset$, then $S_{\{v_1,v_2,v_3\}}=\emptyset$ from Cases 1 and 2.

\vskip 2mm
\noindent
{\bf Claim 3.} \ If $S_{\{v_1,v_2\}}\neq \emptyset$, then $S_{\{v_2,v_3\}}=\emptyset$.

Suppose for a contradiction that $S_{\{v_2,v_3\}}\neq \emptyset$ when $S_{\{v_1,v_2\}}\neq \emptyset$. Then the vertices of $S_{\{v_1,v_2\}}$ (resp., $S_{\{v_2,v_3\}}$) induce a clique, otherwise $\nu(G)\geq 3$, a contradiction. Further, each vertex of $S_{\{v_1,v_2\}}$ is not adjacent to any of $S_{\{v_2,v_3\}}$. If not, one can easily obtain an induced $P_4$, a contradiction. Let $|S_{\{v_1,v_2\}}=a|$ and $|S_{\{v_2,v_3\}}=b|$, then $G$ is isomorphic to $H_5$ in Fig. 2.

\vskip 1.5mm
\noindent
$\bullet$\  Assume that $a=b=\frac{n-3}{2}$, then $d_{v_1}=d_{v_3}=\frac{n-1}{2}$ and $d_{v_2}=n-1$. According to the partition $V(G)=\{V(K_{a+1}), v_2, V(K_{b+1})\}$, $\mathcal{L}(G)$ has an equitable partition and the corresponding quotient matrix is
$$Q_4=
\left(
  \begin{array}{ccc}
    \frac{2}{n-1} & \frac{-\sqrt{2}}{n-1} & 0 \\
    -\frac{1}{\sqrt{2}} & 1 & -\frac{1}{\sqrt{2}} \\
    0 & \frac{-\sqrt{2}}{n-1} & \frac{2}{n-1} \\
  \end{array}
\right).$$
From calculation, the eigenvalues of $Q_4$ are $\{0, \frac{2}{n-1}, \frac{n+1}{n-1}\}$, which are also the eigenvalues of $\mathcal{L}(G)$. Moreover, the multiplicity of $1+\frac{1}{d_{v_1}}=\frac{n+1}{n-1}$ is at least $n-3$ from Lemma \ref{cliquelemma}. Therefore, applying Lemma \ref{quotientmatrixlemma}, we derive that the last unknown $\mathcal{L}$-eigenvalue is
$$n-(n-3)\frac{n+1}{n-1}-\frac{2}{n-1}=\frac{n+1}{n-1}.$$
As a result, the multiplicity of $\frac{n+1}{n-1}$ is $n-2$, contradicting with $G\in \mathcal{G}_1(n, n-3)$.

\vskip 1.5mm
\noindent
$\bullet$\  Assume that $a\neq b$ and $a<b$ without loss of generality. Then we say that $a\leq 2$, otherwise $a\geq 3$, $b\geq 4$ and $G$ contains two distinct  $\mathcal{L}$-eigenvalue with multiplicity at least 3 from Lemma \ref{cliquelemma}, a contradiction. If $a=2$, then the multiplicity of $1+\frac{1}{d_{v_1}}$ is at least 2 by Lemma \ref{cliquelemma}. Since $b=n-5>2$, then $\theta= 1+\frac{1}{d_{v_3}}$ clearly, which yields that the multiplicity of $1+\frac{1}{d_{v_1}}$ is 2. It follows from the trace of $\mathcal{L}(G)$  that
\begin{equation*}
  \begin{array}{rcl}
  n &=& 2(1+\frac{1}{d_{v_1}})+(n-3)1+\frac{1}{d_{v_3}}\\
    &=& 2(1+\frac{1}{3})+(n-3)(1+\frac{1}{n-4})\\
    &=& n+\frac{2}{3}+\frac{1}{n-4}
  \end{array}
\end{equation*}
a contradiction. If $a=1$, with respect to the partition $V(G)=\{V(K_{a+1}), v_2, V(K_{b+1})\}$, then $\mathcal{L}(G)$ has an equitable partition and the corresponding  quotient matrix is
$$Q_5=
\left(
  \begin{array}{ccc}
    1-\frac{1}{2} & \frac{-1}{\sqrt{2(n-1)}} & 0 \\
    \frac{-2}{\sqrt{2(n-1)}} & 1 & \frac{3-n}{(n-1)(n-3)} \\
    0 & \frac{-1}{(n-1)(n-3)} & 1-\frac{n-4}{n-3} \\
  \end{array}
\right),
$$
By computing, the eigenvalues of $Q_5$ are
$$\{0,\ \frac{3n-7\pm 2\sqrt{\frac{n^3+13n^2-125*n+239}{4(n-1)}}}{4(n-3)}\},$$
which are also $\mathcal{L}$-eigenvalues of $G$ from Lemma \ref{quotientmatrixlemma}. Noting that $n\geq 6$ in this case, then $1+\frac{1}{d_{v_1}}=\frac{3}{2}$ and $1+\frac{1}{d_{v_3}}=\frac{n-2}{n-3}$ are two distinct $\mathcal{L}$-eigenvalues of $G$ from Lemma \ref{cliquelemma}. Thus, one can observe that $G$ contains 5 distinct $\mathcal{L}$-eigenvalues, a contradiction.

Combining above discussion, we conclude that $S_{\{v_2,v_3\}}=\emptyset$ when $S_{\{v_1,v_2\}}\neq \emptyset$.

All the proofs are completed. \hfill$\square$

\begin{theorem}\label{theorem2} \ Let $G$ be a connected graph of order $n\geq 5$. Then $G\in \mathcal{G}_1(n, n-3)$ and $G$ is a cograph if and only if $G$ is the graph $G_2$ in Fig. 1.
\end{theorem}

\noindent
{\bf Proof.} \ If $G=G_2$, then clearly $G\in \mathcal{G}_1(n, n-3)$ from Lemma \ref{spectrum} and $G$ is a cograph.

Now suppose that $G\in \mathcal{G}_1(n, n-3)$ and $G$ is a cograph. Also, let $\theta$ be the $\mathcal{L}$-eigenvalue of multiplicity $n-3$, then $\theta\neq 1$ from Lemma \ref{Tian}. In the following, we will show that $G$ must be $G_2$.
Clearly, $G$ cannot be the complete graph $K_n$ and then the diameter $diam(G)=2$.
Let $P_3=v_1v_2v_3$ be a diametrical path of $G$ and for a subset $U$ of $V(P_3)$,
$$S_U=\{u\in V(G)\setminus V(P_3): N_G(u)\cap V(P_3)=U\}.$$
From Lemma \ref{lemma1}, we see that $S_{\{v_1\}}=S_{\{v_2\}}=S_{\{v_3\}}=\emptyset$, and   suppose that $S_{\{v_1,v_2\}}\neq \emptyset$, then $S_{\{v_1,v_3\}}=S_{\{v_2,v_3\}}=S_{\{v_1,v_2,v_3\}}= \emptyset$. Thus all the vertices out of $V(P_3)$ belong to $S_{\{v_1,v_2\}}$. Since $\nu(G)=2$, the vertices of $S_{\{v_1,v_2\}}$ induce a clique of $G$, that is, $G$ is isomorphic to $G_2$.

Next, suppose that $S_{\{v_1,v_2\}}=S_{\{v_2,v_3\}}=\emptyset$. Then all the vertices out of $V(P_3)$ spread in $S_{\{v_1,v_3\}}$ or $S_{\{v_1,v_2,v_3\}}$. In the following, we will show that this cannot hold.

First, assume that $S_{\{v_1,v_3\}}\neq \emptyset$, then we claim that $|S_{\{v_1,v_3\}}|=1$. Otherwise, $|S_{\{v_1,v_3\}}|\geq 2$ and let $v_4, v_5\in S_{\{v_1,v_3\}}$. As $\nu(G)=2$, then $v_4\thicksim v_5$, and thus the vertices $v_i\ (1\leq i\leq 5)$ induce a subgraph isomorphic to $H_1$ in Fig. 2, contradicting with Claim 1 of Lemma \ref{lemma1}. So, $|S_{\{v_1,v_3\}}|=1$ holds and let
$v_4\in |S_{\{v_1,v_3\}}|=1$. Since the order $n\geq 5$ of $G$, then $S_{\{v_1,v_2,v_3\}}\neq \emptyset$. Further, $v_4$ is adjacent to each of $S_{\{v_1,v_2,v_3\}}$. If not, one can also obtain an induced subgraph isomorphic to $H_1$, a contradiction. Now, we can see that $v_2$ and $v_4$ (resp., $v_1$ and $v_3$) are twin points, then the multiplicity of $1$ as an $\mathcal{L}$-eigenvalue is at least 2 from Lemma \ref{twinpointslemma}. As a result, either $\rho_{n-1}(G)= 1$ or $\theta=1$, a contradiction.

Second, assume that $S_{\{v_1,v_3\}}= \emptyset$, then all the vertices of $V(G)\setminus V(P_3)$ belong to $S_{\{v_1,v_2,v_3\}}$. If there are two vertices, say $v_4$ and $v_5$, of $S_{\{v_1,v_2,v_3\}}$ are not adjacent, then we claim that all other vertices of $S_{\{v_1,v_2,v_3\}}$ are adjacent to both $v_4$ and $v_5$. If not, there exists a vertex, say $v_6$, adjacent to exactly one of $v_4$ and $v_5$ (noting that $\nu(G)=2$). As a result, the vertices $v_i\ (2\leq i\leq 6)$ induce an subgraph isomorphic to $H_2$, contradicting with Claim 2 of Lemma \ref{lemma1}. Hence, $v_4$ and $v_5$ are twin points. Noting that $v_1$ and $v_3$ are also twin points, then the multiplicity of $1$ as an $\mathcal{L}$-eigenvalue is at least 2, a contradiction. Consequently, any two of $S_{\{v_1,v_2,v_3\}}$ are adjacent, i.e., $G=K_n-e$. However,  $\rho_{n-1}(K_n-e)= 1$ from \cite{Tian}, a contradiction.

Combining above discussion, we observe that $G$ must be $G_2$.
\hfill$\square$

\begin{remark}\ To complete the characterization of graphs with some normalized Laplacian eigenvalue of multiplicity $n-3$, one need to consider the remaining case, that is, graphs with $\rho_{n-1}(G)\neq 1$ and $\nu(G)= 2$ and $diam(G)=2$. In this case, the edges of graphs are dense. So, it is a challenge to distinguish the edges. However, we conjecture that there is no graphs in this case.
\end{remark}

\vskip 3mm
\noindent
{\large\bf Acknowledgements}\\
The authors thank the anonymous referees for their valuable comments of this paper. This work is supported by the Natural Science Foundation of Shandong Province (No. ZR2019BA016).


{\small
\begin{thebibliography}{90}

\bibitem{vanDam6} E.R. van Dam, G.R. Omidi, Graphs whose normalized Laplacian has three eigenvalues, Linear Algebra Appl. 435 (2011) 2560-2569.
    \bibitem{Guo} J. Li, J.M. Guo, W.C. Shiu, Bounds on normalized Laplacian eigenvalues of graphs, J. Inequal. Appl. 316 (2014) 1-8.
\bibitem{Braga} R.O. Braga, R.R. Del-Vecchio, V.M. Rodrigues, V. Trevisan, Trees with 4 or 5 distinct normalized Laplacian eigenvalues, Linear Algebra Appl. 471 (2015) 615-635.
\bibitem{Guo1} J. Guo, J. Li, W.C. Shiu, The largest normalized Laplacian spectral radius of non-bipartite graphs, Bull. Malaysian Math. Sci. Soc. 39 (1) (2016) 77-87.
\bibitem{Das}  K.C. Das, S. Sun, Normalized Laplacian eigenvalues and energy of trees, Taiwan. J. Math. 20 (3) (2016) 491-507.
\bibitem{Huang2} X. Huang, Q. Huang, On graphs with three or four distinct normalized Laplacian eigenvalues, Algebra Colloquium, 26:1 (2019) 65-82.
\bibitem{Sun1} S. Sun, K.C. Das, On the second largest normalized Laplacian eigenvalue of graphs, Appl. Math. Comput. 348 (2019) 531-541.
\bibitem{Sun2} S. Sun, K.C. Das, Normalized Laplacian spectrum of complete multipartite graphs, Discrete Appl. Math. 284 (2020) 234-245.

\bibitem{Chung} F.R. Chung, Spectral Graph Theory, American Mathematical Society, Providence, RI, 1997.

\bibitem{Fernandes} R. Fernandes, Maria Aguieiras A.de Freitas, Celso M. da Silva Jr., Renata R. Del-Vecchio, Multiplicities of distance Laplacian eigenvalues and forbidden subgraphs, Linear Algebra Appl. 541 (2018) 81-93.
\bibitem{Huang} L. Lu, Q. Huang, X. Huang, On graphs with distance Laplacian spectral radius of multiplicity $n-3$, Linear Algebra Appl. 530 (2017) 485-499.
\bibitem{Tian1} X. Ma, L. Qi, F. Tian, D. wong, Graphs with some distance Laplacian eigenvalue of multiplicity $n-3$, Linear Algebra Appl. 557 (2018) 307-326.
\bibitem{Tian} F. Tian, D. Wong, Characterization of graphs with some normalized Laplacian eigenvalue of multiplicity $n-3$, arXiv:1912.13227.
    
\bibitem{Brouwer} A.E. Brouwer, W.H. Haemers, Spectra of Graphs, Springer, New York, 2012.



\end{thebibliography}
}
\end{document}